\documentclass[12pt, twoside, leqno]{article}
\usepackage{amsmath,amsthm}
\usepackage{amssymb,latexsym}
\usepackage{enumerate}

\newtheorem{thm}{Theorem}[section]

\newtheorem{cor}[thm]{Corollary}
\newtheorem{conj}[thm]{Conjecture}
\newtheorem{lem}[thm]{Lemma}
\newtheorem{prob}[thm]{Problem}

\theoremstyle{definition}

\newtheorem{rem}[thm]{Remark}

\numberwithin{equation}{section}

\usepackage{mathrsfs}
\usepackage{CJK}
\usepackage{amsmath}
\usepackage{fancyhdr}

\begin{document}

\baselineskip=17pt

\title{On simple derivations and the group of polynomial automorphisms commuting with certain derivations}
\author{Dan Yan \footnote{ The author is supported by the NSF of China (Grant No.11601146; 11871241), the China Scholarship Council and the Construct Program of the Key Discipline in Hunan Province.}\\
MOE-LCSM,\\ School of Mathematics and Statistics,\\
 Hunan Normal University, Changsha 410081, China \\
\emph{E-mail:} yan-dan-hi@163.com \\
}
\date{}

\maketitle

\renewcommand{\thefootnote}{}

\renewcommand{\thefootnote}{\arabic{footnote}}
\setcounter{footnote}{0}

\begin{abstract} In the paper, we first study the subgroup of $
K$-automorphisms of $K[x_1,\allowbreak \ldots,x_n]$ which commutes with a simple derivation
of $K[x_1,\ldots,x_n]$. We show that the subgroup of $
K$-automorphisms of $K[x_1,\ldots,x_n]$ which commutes with a simple derivation
of $K[x_1,\ldots,x_n]$ consists of translations under certain hypothesis. We also prove that the
subgroup of $K$-automorphisms of $K[x_1,x_2]$ which commutes with
a simple derivation is trivial under the same hypothesis, and the subgroup of $ K$-automorphisms of $K[x_1,\ldots,x_n]$ which commutes with simple Shamsuddin derivations is trivial under the same hypothesis. Then we give an affirmative answer to the conjecture proposed in \cite{13} for $n=2$.
\end{abstract}
\section{Introduction}

Throughout this paper, we will write $K$ for any field with characteristic zero and $K[x]=K[x_1,\allowbreak\ldots,x_n]$ for the polynomial algebra over $K$
in $n$ indeterminates $x_1,x_2,\ldots,x_n$. $\partial_{x_i}$ will
denote the derivation $\frac{\partial}{\partial x_i}$ of $ K[x]$ for
all $1\leq i\leq n$. We abbreviate $\frac{\partial f_j}{\partial x_i}$ as $f_{jx_i}$.

A $K$-derivation $D:K[x]\rightarrow K[x]$ of $ K[x]$ is a $
K$-linear map such that
$$D(ab)=D(a)b+aD(b)$$
for any $a,b\in  K[x]$, and $D(c)=0$ for any $c\in K$. The set of all $ K$-derivations of $K[x]$ is denoted by $\operatorname{Der}_ K( K[x])$. An ideal $I$ of $
K[x]$ is called $D$-$stable$ if $D(I)\subset I$. $ K[x]$ is called
$D$-$simple$ if it has no proper nonzero $D$-$stable$ ideal. The $
K$-derivation $D$ is called $simple$ if $ K[x]$ has no $D$-$stable$
ideals other than $0$ and $ K[x]$. For some examples of simple
derivations, see \cite{1}, \cite{4}, \cite{5}, \cite{6}, \cite{7},
\cite{10}, \cite{8}.

A polynomial map is a map $f=(f_1,f_2,\ldots,f_n): K^n\rightarrow
 K^n$ of the form
$$(x_1,x_2,\ldots,x_n)\rightarrow (f_1(x_1,\ldots,x_n),\ldots,f_n(x_1,\ldots,x_n)),$$
where each $f_i\in K[x]$. Such a polynomial map is called invertible
if there exists a polynomial map $h=(h_1,h_2,\ldots,h_n):
K^n\rightarrow K^n$ such that $x_i=h_i(f_1,\ldots,f_n)$ for all
$1\leq i\leq n$. Invertible polynomial maps correspond one-to-one
with $ K$-automorphisms of the polynomial ring $ K[x]$ given by
$$f\rightarrow f^*:g\rightarrow g(f).$$
So describing invertible polynomial maps from $K^n\rightarrow K^n$
is the same as describing $ K$-automorphisms of $ K[x]$. The group
of $ K$-automorphisms of $ K[x]$ will be denoted by
$\operatorname{Aut}( K[x])$. For more information about invertible
polynomial maps, see \cite{12}.

 Let
$\operatorname{Aut}( K[x])$ act on $\operatorname{Der}_ K( K[x])$
as:
$$(\rho,D)\rightarrow\rho^{-1}\circ D\circ\rho=\rho^{-1} D\rho.$$
The isotropy subgroup is defined to be:
$$\operatorname{Aut}(
K[x])_D:=\{\rho\in\operatorname{Aut}( K[x])|\rho D=D\rho\}.$$

In \cite{2} and \cite{3}, the authors proposed the following
conjecture.

\begin{conj}
If $D$ is a simple derivation of $ K[x_1,x_2]$, then
$\operatorname{Aut}( K[x_1,x_2])_D$ is finite.
\end{conj}

We can ask the same question for the polynomial ring in $n$
variables $ K[x]$. However, the same question is not true for other
cases. If $n=1$, then all the simple derivations of $ K[x_1]$ have
the form $c\cdot\partial_1$ for any $c\in  K^*$. Thus, it is easy to
compute that $\operatorname{Aut}( K[x_1])_D=\{x_1+c_1|
~c_1\in  K$\}. Therefore, $\operatorname{Aut}(K[x_1])_D$ is infinite.
If $n\geq 3$, then let
$D=(1-x_1x_2)\partial_1+x_1^3\partial_2+x_2\partial_3+\cdots+x_{n-1}\partial_n$,
it follows from \cite{9} that $D$ is simple. It is easy to check
that
$(x_1,x_2,\ldots,x_{n-1},x_n+c_n)\in \operatorname{Aut}( K[x])_D$
for any $c_n\in  K$. Thus, $\operatorname{Aut}(
K[x])_D$ is infinite.\\

Although the conjecture is not true for $n$ variables with $n\geq
3$, we can still ask the same question for simple Shamsuddin
derivations in $n$ variables, where $n\geq 3$.

A derivation $D$ of $ K[x]$ is said to be a Shamsuddin derivation if
$D=\partial_{x_1}+\sum_{i=2}^n(a_ix_i+b_i)\partial_{x_i}$ with $a_i,
b_i\in  K[x_1]$ for all $2\leq i\leq n$.

\begin{conj}
 If $D$ is a simple Shamsuddin derivation
of $ K[x]$, then $\operatorname{Aut}( K[x])_D\allowbreak =\{id\}$.
\end{conj}

Rene Baltazar shows that Conjecture 1.1 is true for the Shamsuddin
derivations of $ K[x_1,x_2]$ in \cite{3}. In particular, the author
proves that $\operatorname{Aut}( K[x_1,x_2])_D=\{id\}$ if $D$ is a
simple Shamsuddin derivation of $ K[x_1,x_2]$. Thus, conjecture 1.2
has an affirmative answer if $n=2$. L.G.Mendes and Ivan Pan have proved that the subgroup of $K$-automorphisms of $ K[x_1,x_2]$ which commute with simple
derivations is trivial, see \cite{11}. In \cite{13}, L.N.Bertoncello and
D.Levcovitz also proposed Conjecture 1.2 and have proved that the subgroup of $ K$-automorphisms of $K[x]$ which commute with simple Shamsuddin derivations is trivial.
They also proposed the following conjecture.

\begin{conj}
 If $D$ is a Shamsuddin derivation
of $K[x]$, then $D$ is simple if, and only if $\operatorname{Aut}( K[x])_D\allowbreak =\{id\}$.
\end{conj}

In our paper, we give a different point of view to the first two conjectures and pose a conjecture for all simple derivations of $K[x]$. In addition, we give an affirmative answer to the Conjecture 1.3 for $n=2$.

We arrange as follow: we prove that $\operatorname{Aut}(
K[x])_D=\{x+c|~c\in K^n$\} if $D$ is simple and $f-f(0)\in \operatorname{Aut}(
K[x])_D$ in the case that $f\in \operatorname{Aut}(
K[x])_D$. Then, in section 3, we show that
$\operatorname{Aut}(K[x])_D=\{id\}$ if $D$ is a simple derivation of
the form:
$D=p(x_1)\partial_{x_1}+\sum_{i=2}^nq_i(x_1,x_i)\partial_{x_i}$ with
$p(x_1)\in K[x_1]$, $q_i(x_1,x_i)\in  K[x_1,x_i]$ for all $2\leq i\leq n$ under the hypothesis that $f-f(0)\in \operatorname{Aut}(
K[x])_D$ if $f\in \operatorname{Aut}(K[x])_D$. We also give an affirmative answer to the conjecture 1.3 for $n=2$ in section 3.

\section{Affirmative answer to conjecture 1.1 under certain hypothesis}

In the section, we show that $\operatorname{Aut}(
K[x_1,x_2])_D=\{id\}$ if $D$ is a simple derivation over $
K[x_1,x_2]$ and $f-f(0)\in \operatorname{Aut}(
K[x_1,x_2])_D$ in the case that $f\in \operatorname{Aut}(
K[x_1,x_2])_D$.

\begin{thm}
Let $D=\sum_{i=1}^np_i(x)\partial_{x_i}$ be a simple derivation over $
K[x]$, where $p_i(x)\in K[x]$ for all $1\leq i\leq n$. If $f$, $f-f(0)\in \operatorname{Aut}(K[x])_D$ then $f=x+f(0)$. In particular, $f$ is a translation.
\end{thm}
\begin{proof}
If $f\in \operatorname{Aut}(
K[x])_D$, then $f=(f_1,f_2,\ldots,f_n)$ is an invertible polynomial map. Let $\bar{f}=f-f(0)$. Then $\bar{f}$ is an invertible polynomial map with $\bar{f}(0)=0$. Thus, we can assume that
$h=(h_1,\ldots,h_n)$ be the inverse of $\bar{f}$. By the definition of
invertible polynomial map, we have $x_i=h_i(\bar{f}_1,\ldots,\bar{f}_n)$ for all
$1\leq i\leq n$. Substituting $x_1=\cdots=x_n=0$ to the above
equation, we have
$0=h_i(\bar{f}_1(0,\ldots,0),\ldots,\bar{f}_n(0,\ldots,0))=h_i(0,\ldots,0)$ for
all $1\leq i\leq n$. Thus, we can write
$x_i=\sum_{j=1}^nh_i^{(j)}(\bar{f}_1,\ldots,\bar{f}_n)\cdot \bar{f}_j$ for all $1\leq
i\leq n$, where $h_i^{(j)}(\bar{f}_1,\ldots,\bar{f}_n)\in  K[x]$. Therefore, the
ideal $(x_1,\dots,x_n)$ is equal to the ideal $(\bar{f}_1,\ldots,\bar{f}_n)$.
Let ${\bf m}$ be the ideal $(x_1,\dots,x_n)$ and $\rho(x_i)=\bar{f}_i(x_1,\ldots,x_n)$ for all $1\leq i\leq n$.
Then we have $\rho({\bf m})={\bf m}$. Since $\rho \in \operatorname{Aut}(
K[x])_D$, it follows from Proposition 7
in \cite{3} that $\rho=id$. That is, $f=x+f(0)$, where $f(0)=(f_1(0),f_2(0),\ldots,f_n(0))$.
\end{proof}

\begin{cor}
Let $D$ be a simple derivation of $K[x]$ and $S=\{f\in \operatorname{Aut}(K[x])_D|\allowbreak f(0)=0\}$. Then $S=\{id\}$.
\end{cor}
\begin{proof}
Suppose $f\in S$. Then $f-f(0)\in S$. Thus, the conclusion follows from Theorem 2.1.
\end{proof}

\begin{lem}
Let $D=\sum_{i=1}^np_i(x)\partial_{x_i}$ be a derivation of $
K[x]$, where $p_i(x)\in K[x]$ for all $1\leq i\leq n$. Then the following two statements are equivalent:

$(1)$ $f-f(0)\in \operatorname{Aut}(K[x])_D$ in the case that $f\in \operatorname{Aut}(K[x])_D$;

$(2)$ $p_i(x_1+c_1,\ldots,x_n+c_n)=p_i(x_1,\ldots,x_n)$ for $c_j=-f_j(0) \in
K$ and for all $1\leq i,j\leq n$.
\end{lem}
\begin{proof}
Suppose that $(1)$ is true and $f \in \operatorname{Aut}(K[x])_D$. Then $f-f(0)\in \operatorname{Aut}(K[x])_D$. Thus, we have the following equations:
$$\sum_{i=1}^np_i(x)f_{jx_i}=p_j(f)$$
and
$$\sum_{i=1}^np_i(x)f_{jx_i}=p_j(f-f(0))$$
for all $1\leq j\leq n$.
Therefore, we have $$p_j(f)=p_j(f-f(0))$$
by comparing the above two equations for all $1\leq j\leq n$. Composing the inverse of $f$ from the right to the two sides of the above equation, we have $p_j(x)=p_j(x-f(0))$ for all $1\leq j\leq n$. Let $c=-f(0)$. Then the statement $(2)$ follows.

Conversely, since $c=-f(0)$, we have $$p_j(x)=p_j(x-f(0))$$
 for all $1\leq j\leq n$. Composing $f$ from the right to the two sides of the above equation, we have $p_j(f)=p_j(f-f(0))$. If $f\in \operatorname{Aut}(K[x])_D$, then we have $\sum_{i=1}^np_i(x)f_{jx_i}=p_j(f)$. Thus, we have $\sum_{i=1}^np_i(x)f_{jx_i}=p_j(f-f(0))$ by combining the former two equations. That is, $f-f(0)\in \operatorname{Aut}(K[x])_D$.
\end{proof}

\begin{lem}
Let $p(x)$ be a polynomial over $K[x]$ and
$p(x_1+c_1,\ldots,x_n+c_n)=p(x_1,\ldots,x_n)$ for some $c_i\in K$ and for all
$1\leq i\leq n$. If there exists $c_k\neq 0$ for some $k\in
\{1,2,\ldots,n\}$, then $p(x)\in
K[c_kx_1-c_1x_k,\ldots,c_kx_{k-1}-c_{k-1}x_k,c_kx_{k+1}-c_{k+1}x_k,\ldots,\allowbreak c_kx_n-c_nx_k]$.
\end{lem}
\begin{proof}
Let $p^{(m)}$ be the highest homogeneous part of $p$ with $\deg p^{(m)}=m$. Since
$$p(x_1+c_1,\ldots,x_n+c_n)=p(x_1,\ldots,x_n),~~~(2.1)$$
we have
$$\sum_{i=1}^nc_ip_{x_i}^{(m)}=0~~~~~~~~~~~~~~~~~~~(2.2)$$
by comparing the homogeneous part of degree $m-1$ of equation
$(2.1)$. Let $\bar{x}_j=c_kx_j-c_jx_k$ for all $1\leq j\leq n$,
$j\neq k$; $\bar{x}_k=x_k$. Thus, we have
$$p_{x_j}^{(m)}=c_kp_{\bar{x}_j}^{(m)}~~~~~~~~~~~~~~~~~~~(2.3)$$
for all $1\leq j\leq n$, $j\neq k$ and
$$p_{x_k}^{(m)}=(c_k+1)p_{\bar{x}_k}^{(m)}-\sum_{i=1}^nc_ip_{\bar{x}_i}^{(m)}.~~~~(2.4)$$
Thus, we have $p_{\bar{x}_k}^{(m)}=0$ by substituting equations $(2.3)$, $(2.4)$ to equation $(2.2)$.  That is, $p^{(m)}\in
K[\bar{x}_1,\ldots,\bar{x}_{k-1},\bar{x}_{k+1},\ldots,\bar{x}_n]$.
Clearly, we have
$p^{(m)}(x_1+c_1,\ldots,x_n+c_n)=p^{(m)}(x_1,\ldots,x_n)$. Since
$p=p^{(m)}+p^{(m-1)}+\cdots+p^{(1)}+p^{(0)}$ and
$\bar{p}(x_1+c_1,\ldots,x_n+c_n)=\bar{p}(x_1,\ldots,x_n)$, where
$\bar{p}=p-p^{(m)}$ and $p^{(l)}$ is the homogeneous part of $p$
with $\deg p=l$, $0\leq l\leq m$, we can continue the process for
$\bar{p}$. Then we have $p^{(l)}\in
K[\bar{x}_1,\ldots,\bar{x}_{k-1},\bar{x}_{k+1},\ldots,\bar{x}_n]$
for all $0\leq l\leq m-1$. That is, $p(x)\in
K[c_kx_1-c_1x_k,\ldots,c_kx_{k-1}-c_{k-1}x_k,c_kx_{k+1}-c_{k+1}x_k,\ldots, c_kx_n-c_nx_k]$.
\end{proof}

\begin{rem}
Wenhua Zhao point out that Lemma 2.4 may be already known, but we have not found the reference. Thus, we give our proof here. He also gives another proof of Lemma 2.4.
\end{rem}
Another proof of Lemma 2.4.

Without loss of generality, we can assume that $k=1$. Let $\bar{x}_j=c_1x_j-c_jx_1$ for all $1\leq j\leq n$,
$j\neq 1$; $\bar{x}_1=x_1$. Thus, we have $p(x_1,\ldots,x_n)=\bar{p}(\bar{x}_1,\bar{x}_2,\ldots,\bar{x}_n)$ for some $\bar{p}(\bar{x}_1,\bar{x}_2,\ldots,\bar{x}_n)\in K[\bar{x}_1,\ldots,\bar{x}_n]$. Therefore, we have $p(x_1+c_1,\ldots,x_n+c_n)=\bar{p}(\bar{x}_1+c_1,\bar{x}_2\ldots,\bar{x}_n)$. Then we have $$\bar{p}(\bar{x}_1,\bar{x}_2,\ldots,\bar{x}_n)=\bar{p}(\bar{x}_1+c_1,\bar{x}_2\ldots,\bar{x}_n).~~~(2.5)$$
for $c_1\neq 0$. We can view $\bar{p}$ as a polynomial in $R[\bar{x}_1]$ with coefficients in $R$, where $R=K[\bar{x}_2,\ldots,\bar{x}_n]$. Thus, we have $\bar{p}\in K[\bar{x}_2,\ldots,\bar{x}_n]$ by comparing the degree of $\bar{x}_1$ of equation $(2.5)$. Then the conclusion follows.

\begin{lem}
Let $D=\sum_{i=1}^np_i(x)\partial_{x_i}$ be a simple derivation of
$ K[x]$, where $p_i(x)\in K[x]$ for all $1\leq i\leq n$. Then we
have the following statements:

$(1)$ $p_1, p_2, \ldots, p_n$ are linearly independent;

$(2)$ $p_i(x)\notin  K[x_i]$ or $p_i(x)\in  K^*$ for all $1\leq
i\leq n$;

$(3)$ $p_i(x_1,\ldots,x_{i-1},0,x_{i+1},\ldots,x_n)\neq 0$ for all $1\leq i\leq n$.
\end{lem}
\begin{proof}
$(1)$ Suppose that $p_1, p_2, \ldots, p_n$ are linearly dependent.
Then there exists
$(\lambda_1,\lambda_2,\allowbreak\ldots,\lambda_n)\in K^n$ and
$(\lambda_1,\lambda_2,\ldots,\lambda_n)\neq (0,0,\ldots,0)$ such
that $\sum_{i=1}^n\lambda_ip_i(x)=0$. Thus, the ideal
$(\sum_{i=1}^n\lambda_ix_i)$ is $D$-stable. This is a contradiction.
Thus, $p_1, p_2, \ldots,\allowbreak p_n$ are linearly independent.

$(2)$ If $p_i(x)\in  K[x_i]$ and $p_i(x_1,\ldots,x_n)\notin K^*$,
then the ideal $(p_i(x))$ is
$D$-stable. This is a contradiction. Thus, $p_i(x)\notin K[x_i]$ or
$p_i(x)\in  K^*$ for all $1\leq i\leq n$.

$(3)$ If $p_i(x_1,\ldots,x_{i-1},0,x_{i+1},\ldots,x_n)=0$, then the ideal $(x_i)$ is $D$-stable. This is a contradiction. Thus, $p_i(x_1,\ldots,x_{i-1},0,x_{i+1},\ldots,x_n)\neq 0$ for all $1\leq i\leq n$.
\end{proof}

\begin{thm}
Let $D=p_1(x_1,x_2)\partial_{x_1}+p_2(x_1,x_2)\partial_{x_2}$ be a derivation of $K[x_1,x_2]$. If $p_i(x+c)=p_i(x)$ for some $c\in K^2 \setminus \{(0,0)\}$ and for all $1\leq i\leq 2$, then $D$ is not simple.
\end{thm}
\begin{proof}
Suppose that $D$ is a simple derivation and, without loss of generality, $c_1\neq 0$. Then it follows from Lemma 2.4 that $p_1, p_2\in K[c_1x_2-c_2x_1]$. Let
$\bar{x}_1=x_1$, $\bar{x}_2=c_1x_2-c_2x_1$. Then
$\partial_{x_1}=\partial_{\bar{x}_1}-c_2\partial_{\bar{x}_2}$ and
$\partial_{x_2}=c_1\partial_{\bar{x}_2}$. Thus, we have
$$D=r_1(\bar{x}_2)\partial_{\bar{x}_1}+(c_1r_2(\bar{x}_2)-c_2r_1(\bar{x}_2))\partial_{\bar{x}_2},$$
where $r_1(\bar{x}_2)=p_1(\bar{x}_2), r_2(\bar{x}_2)=p_2(\bar{x}_2)$, which are polynomials in $K[\bar{x}_2]$.
Since $ K[x_1,x_2]= K[\bar{x}_1,\bar{x}_2]$, $D$ is a simple
derivation over $K[\bar{x}_1,\bar{x}_2]$. It follows from Lemma 2.6
(2) that $c_1r_2(\bar{x}_2)-c_2r_1(\bar{x}_2)\in  K^*$. Let
$\bar{c}=c_1r_2(\bar{x}_2)-c_2r_1(\bar{x}_2)$ and $h(\bar{x}_2)\in
 K[\bar{x}_2]$ such that $h'(\bar{x}_2)=r_1(\bar{x}_2)$. Then
the ideal $(\bar{x}_1-\bar{c}^{-1}h(\bar{x}_2))$ is a $D$-$stable$
ideal of $ K[\bar{x}_1,\bar{x}_2]$. This is a contradiction. Thus,
$D$ is not simple.
\end{proof}

\begin{cor}
Let $D=p_1(x_1,x_2)\partial_{x_1}+p_2(x_1,x_2)\partial_{x_2}$ be a simple derivation of $K[x_1,x_2]$. If $f-f(0)\in \operatorname{Aut}(
K[x_1,x_2])_D$ in the case that $f\in \operatorname{Aut}(
K[x_1,x_2])_D$, then $\operatorname{Aut}(
K[x_1,x_2])_D=\{id\}$.
\end{cor}
\begin{proof}
Let $\rho \in \operatorname{Aut}(K[x_1,x_2])_D$ and $\rho(x_1)=f_1(x_1,x_2)$, $\rho(x_2)=f_2(x_1,x_2)$. It follows from Lemma 2.3 that $p_i(x_1-f_1(0),x_2-f_2(0))=p_i(x_1,x_2)$ for all $1\leq i\leq 2$. Since $D$ is simple, it follows from Theorem 2.7 that $f_1(0)=f_2(0)=0$. Then it follows from Theorem 2.1 and Corollary 2.2 that $\rho=id$. That is, $\operatorname{Aut}(K[x_1,x_2])_D=\{id\}$.
\end{proof}

\section{Affirmative answer to Conjecture 1.3 for n=2 and conjecture 1.2 under certain hypothesis}

In the section, we show that $\operatorname{Aut}(K[x])_D=\{id\}$ if
$D$ is simple and of the form
$p(x_1)\partial_{x_1}\allowbreak+\sum_{i=2}^nq_i(x_1,x_i)\partial_{x_i}$,
where $p(x_1) \in K[x_1]$, $q_i(x_1,x_i)\in K[x_1,x_i]$ for all $2\leq i\leq n$, under the hypothesis that $f-f(0)\in \operatorname{Aut}(K[x])_D$ if $f\in \operatorname{Aut}(K[x])_D$. In addition, we give an affirmative answer to conjecture 1.3 for $n=2$.

\begin{thm}
Let $D=p(x_1)\partial_{x_1}+\sum_{i=2}^nq_i(x_1,x_i)\partial_{x_i}$ be
a derivation of $K[x]$, where $n\geq 2$  and $p(x_1) \in K[x_1]$, $q_i(x_1,x_i)\in K[x_1,x_i]$ for all $2\leq i\leq n$. If $q_i(x_1+c_1,x_i+c_i)=q_i(x_1,x_i)$ for some $c=(c_1,c_2,\ldots,c_n)\in K^n\setminus \{(0,0,\ldots,0)\}$ and for all $2\leq i\leq n$, then $D$ is not simple.
\end{thm}
\begin{proof}
If $p(x_1)=0$, then the ideal $(x_1)$ is $D$-stable. Thus, $D$ is not simple.

If $p(x_1)\in K[x_1]\setminus K$, then the ideal $(p(x_1))$ is $D$-stable. Thus, $D$ is not simple.

If $p(x_1)\in K^*$, then let $e:=p(x_1)$. Since
$$q_i(x_1+c_1,x_i+c_i)=q_i(x_1,x_i)~~~~~~~~~~~~~~~~~~~~~(3.1)$$
for some $(c_1,c_2,\ldots,c_n)\in K^n\setminus \{(0,0,\ldots,0)\}$ and for all $2\leq i\leq n$, we have the following two cases.

(1) If $c_1\neq 0$, then it follows from Lemma 2.4 that $q_i\in K[c_1x_i-c_ix_1]$ for all $2\leq i\leq n$. Let $\bar{x}_1=x_1$, $\bar{x}_i=c_1x_i-c_ix_1$ for all $2\leq i\leq n$. Then $\partial_{x_1}=\partial_{\bar{x}_1}-\sum_{i=2}^nc_i\partial_{\bar{x}_i}$, $\partial_{x_i}=c_1\partial_{\bar{x}_i}$ for all $2\leq i\leq n$. Thus, we have
$$D=e\partial_{\bar{x}_1}+\sum_{i=2}^n(c_1q_i(\bar{x}_i)-ec_i)\partial_{\bar{x}_i}.$$

If there exists $i_0\in \{2,3,\ldots,n\}$ such that $\deg q_{i_0}(\bar{x}_{i_0})\geq 1$, then the ideal $(c_1q_{i_0}(\bar{x}_{i_0})-ec_{i_0})$ is $D$-stable. Thus, $D$ is not simple.

If $q_i(\bar{x}_i)\in K$ for all $2\leq i\leq n$, then let $q_i:=q_i(\bar{x}_i)$, then the ideal $((c_1q_i-ec_i)\bar{x}_1-e\bar{x}_i)$ is $D$-stable. Thus, $D$ is not simple.

(2) If $c_1=0$, then equation (3.1) has the following form:
$$q_i(x_1,x_i+c_i)=q_i(x_1,x_i)~~~~~~~~~~~~~~~~~~~~~(3.2)$$
for some $(c_1,c_2,\ldots,c_n)\in K^n\setminus \{(0,0,\ldots,0)\}$ and for all $2\leq i\leq n$.
Suppose that $\deg_{x_i}q_i(x_1,x_i)=m_i$ for $2\leq i\leq n$.

If $m_i\geq 1$ for all $2\leq i\leq n$, then we have $c_i=0$ by comparing the coefficients of $x_i^{m_i-1}$ of equation $(3.2)$ for all $2\leq i\leq n$:contradiction.

If there exists $\tilde{i}\in \{2,3,\ldots,n\}$ such that $m_{\tilde{i}}=0$, then we have $q_{\tilde{i}}(x_1,x_{\tilde{i}})\in K[x_1]$. Let $q_{\tilde{i}}(x_1):=q_{\tilde{i}}(x_1,x_{\tilde{i}})$ and $Q(x_1)=\int q_{\tilde{i}}(x_1)dx_1$. Then the ideal $(ex_{\tilde{i}}-Q(x_1))$ is $D$-stable. Thus, $D$ is not simple. This completes the proof.
\end{proof}

\begin{cor}
Let $D=p(x_1)\partial_{x_1}+\sum_{i=2}^nq_i(x_1,x_i)\partial_{x_i}$ be
a simple derivation of $K[x]$, where $n\geq 2$ and $p(x_1)\in K[x_1]$, $q_i(x_1,x_i)\in K[x_1,x_i]$ for all $2\leq i\leq n$. If $f-f(0)\in \operatorname{Aut}(K[x])_D$ in the case that $f\in \operatorname{Aut}(K[x])_D$, then $\operatorname{Aut}(K[x])_D=\{id\}$.
\end{cor}
\begin{proof}
Let $\rho\in \operatorname{Aut}(K[x])_D$ and $\rho(x_i)=f_i$ for all $1\leq i\leq n$. It follows from Lemma 2.3 that $p(x_1-f_1(0))=p(x_1)$ and $q_i(x_1-f_1(0),x_i-f_i(0))=q_i(x_1,x_i)$ for all $2\leq i\leq n$. Since $D$ is simple, it follows from Theorem 3.1 that $f_1(0)=\cdots=f_n(0)=0$. Then it follows from Theorem 2.1 and Corollary 2.2 that $\rho=id$. That is, $\operatorname{Aut}(K[x])_D=\{id\}$.
\end{proof}

\begin{cor}
Let
$D=\partial_{x_1}+\sum_{i=2}^n(a_i(x_1)x_i+b_i(x_1))\partial_{x_i}$
be a simple Shamsuddin derivation, where $n\geq 2$ and $a_i(x_1), b_i(x_1) \in
K[x_1]$ for all $2\leq i\leq n$. If $f-f(0)\in \operatorname{Aut}(K[x])_D$ in the case that $f\in \operatorname{Aut}(K[x])_D$, then $\operatorname{Aut}(K[x])_D=\{id\}$.
\end{cor}
\begin{proof}
The conclusion follows from Corollary 3.2.
\end{proof}

We proposed the following conjecture for all simple derivations.

\begin{conj}
Let $D=\sum_{i=1}^np_i(x)\partial_{x_i}$ be a derivation over $
K[x]$, where $p_i(x)\in K[x]$ for all $1\leq i\leq n$. If $D$ is a simple derivation, then $\operatorname{Aut}(K[x])_D=T$, where $T$ is the translation subgroup of $\operatorname{Aut}(K[x])$.
\end{conj}

\begin{rem}
If Conjecture 3.4 has an affirmative answer, then Conjecture 1.1 and Conjecture 1.2 have affirmative answers and the condition in Corollary 2.8 and Corollary 3.2 that $f-f(0)\in \operatorname{Aut}(K[x])_D$ in the case that $f\in \operatorname{Aut}(K[x])_D$ can be removed. In addition, Conjecture 3.4 has an affirmative answer if $n=1$. It follows from Theorem 2.1 that we only need to prove the following statement in order to prove Conjecture 3.4.
\end{rem}

\begin{prob}
Let $D$ be a simple derivation of $K[x]$, If $f\in \operatorname{Aut}(K[x])_D$, then $f-f(0)\in \operatorname{Aut}(K[x])_D$.
\end{prob}

Next we give an affirmative answer to Conjecture 1.3 for $n=2$.

\begin{lem}
Let $D=\partial_{x_1}+b(x_1)\partial_{x_2}$ be a derivation of $K[x_1,x_2]$ with $b(x_1)\in K[x_1]$. Then $\operatorname{Aut}(K[x_1,x_2])_D=\{(f,g)\}$, where $f=x_1+p(x_2-h(x_1))$, $g=\sum_{k=0}^m\frac{1}{k+1}C_kf^{k+1}+\bar{q}(x_2-h(x_1))$ and $b(x_1)=\sum_{k=0}^mC_kx_1^k$, $h(x_1)=\int b(x_1)dx_1$, $p(x_2-h(x_1))\in K[x_2-h(x_1)]$, $\bar{q}(x_2-h(x_1))=\tilde{c}\cdot(x_2-h(x_1))+\bar{c}$ for some $\tilde{c}\in K^*$, $\bar{c}\in K$.
\end{lem}
\begin{proof}
Let $\rho \in \operatorname{Aut}(K[x_1,x_2])_D$ with $\rho(x_1)=f(x_1,x_2)$, $\rho(x_2)=g(x_1,x_2)$. Then we have the following equations:
$$f_{x_1}+b(x_1)f_{x_2}=1~~~~~~~~~~~~~~~~~~~~(3.4)$$
and
$$g_{x_1}+b(x_1)g_{x_2}=b(f).~~~~~~~~~~~~~~~~~~~(3.5)$$
Let $\tilde{x}_1=x_1$, $\tilde{x}_2=x_2-\int b(x_1)dx_1$. It follows from equation (3.4) that $f_{\tilde{x}_1}=1$. Thus, we have $f=\tilde{x}_1+p(\tilde{x}_2)$ for some $p(\tilde{x}_2)\in K[\tilde{x}_2]$. That is, $f(x_1,x_2)=x_1+p(x_2-h(x_1))$, where $h(x_1)=\int b(x_1)dx_1$. Since $b(x_1)=\sum_{k=0}^mC_kx_1^k$, it follows from equation (3.5) that $g_{\tilde{x}_1}=\sum_{k=0}^mC_kf^k$. Thus, we have
$$g=\sum_{k=0}^m\frac{1}{k+1}C_kf^{k+1}+\bar{q}(\tilde{x}_2)$$
for some $\bar{q}(\tilde{x}_2)\in K[\tilde{x}_2]$. Since $(f,g)\in \operatorname{Aut}(K[x_1,x_2])$, we have $\det J(f,g)=\tilde{c}\in K^*$. Thus, we have $\bar{q}'(x_2-h(x_1))=\tilde{c}$. Therefore, we have
$$\bar{q}(x_2-h(x_1))=\tilde{c}\cdot(x_2-h(x_1))+\bar{c}$$
for some $\bar{c}\in K$, $\tilde{c}\in K^*$. This completes the proof.
\end{proof}

\begin{rem}
Since $b(x_1+c_1)=b(x_1)$ for any $c=(c_1,c_2)\in K^2$ with $c_1=0$, it follows from Theorem 2.7 that the derivation $D$ in Lemma 3.6 is not simple.
\end{rem}

\begin{thm}
Let $D=\partial_{x_1}+(a(x_1)x_2+b(x_1))\partial_{x_2}$ be a Shamsuddin derivation of $K[x_1,x_2]$. Then the following two statements are equivalent:

(1) $D$ is a simple derivation;

(2) $\operatorname{Aut}(K[x_1,x_2])_D=\{id\}$.
\end{thm}
\begin{proof}
$(1)\Rightarrow (2)$ It follows from Theorem 6 in \cite{3}.

$(2)\Rightarrow (1)$ If $a(x_1)\neq 0$, then it follows from Theorem 2 in \cite{11} that $D$ is a simple derivation. If $a(x_1)=0$, then the conclusion follows from Lemma 3.6.
\end{proof}

\begin{rem}
Since Conjecture 1.3 has an affirmative answer for $n=2$, we want to know whether the converse of Conjecture 3.4 is true or not. In general, it is not true. Let $D=(d_2x_1^2+d_1x_1+d_0)\partial_{x_1}$ be a derivation over $K[x_1]$ with $d_2\neq 0$. It is easy to compute that $\operatorname{Aut}(K[x_1])_D=\{id\}$ and the ideal $(d_2x_1^2+d_1x_1+d_0)$ is $D$-stable. Thus, $D$ is not simple. Let $D=\partial_{x_1}+(a_2(x_1)x_2^2+a_1(x_1)x_2)\partial_{x_2}$ be a derivation over $K[x_1,x_2]$ with $\deg a_2(x_1)\geq 1$. Then the ideal $(x_2)$ is $D$-stable. Thus, $D$ is not simple. It follows from Theorem 4.2 (i) and (ii) in \cite{13} that $\operatorname{Aut}(K[x_1,x_2])_D=\{id\}$. Thus, the above two examples are counterexamples of the converse of Conjecture 3.4 for $n=1,~2$.
\end{rem}

{\bf{Acknowledgement}}: The author is very grateful to professor Wenhua Zhao for pointing out some mistakes and misprints of the paper. She is also grateful to the Department of Mathematics of Illinois State University, where this paper was finished, for hospitality during her stay as a visiting scholar.

\end{document}